\documentclass{article}
\usepackage{amssymb}

\usepackage[english,francais]{babel}

\newtheorem{e-proposition}[theorem]{Proposition}

\newtheorem{e-definition}[theorem]{Definition\rm}

\newtheorem{theoreme}{Th\'eor\`eme}[section]

\newtheorem{proposition}[theoreme]{Proposition}

\def \H {{\mathbb H}}

\font \m=msbm10

\def\H{{\hbox {\m H}}}

\begin{document}
\title{SLEs as boundaries of clusters of Brownian loops}
\author{Wendelin Werner}

\date{Universit\'e Paris-Sud and IUF,
Laboratoire de Math\'ematiques, 91405 Orsay cedex, France}
\maketitle

\begin{abstract}
In this research announcement, we show that SLE curves 
can in fact be viewed as boundaries of certain 
clusters of Brownian loops (of the clusters in a Brownian loop soup).
For small densities $c$ of loops, we
 show that  the outer boundaries of the clusters created by the  Brownian 
loop soup are SLE$_\kappa$-type curves where 
$\kappa \in (8/3, 4]$ and $c$ related by the usual relation $c=(3\kappa-8)(6-\kappa)/2\kappa$
(i.e. $c$ corresponds to the central charge of the model).
This gives (for any Riemann surface) a simple construction of a natural countable family of
random disjoint SLE$_\kappa$ loops,
that behaves ``nicely'' under perturbation of the surface
and is related to various aspects of conformal field theory and representation theory.
\end {abstract}

\section {Background}

The goal of this paper is announce some results that relate 
the Brownian loop soup introduced in \cite {LWls} to the Schramm-Loewner Evolutions (SLE)
(in particular for values of the parameter between $8/3$ and $4$) and random families of disjoint 
SLE loops (as they might appear 
in the scaling limit of many $2d$ statistical physics models, and in conformal field theory). A
more complete paper \cite {Wip} on the same subject (with proofs, that discusses also various consequences of this approach)
 is in preparation.

SLE processes have been introduced by Schramm in \cite {Schramm1}, building on the observation that 
they are the only processes that have a certain (conformal) Markovian-type property. There is a one-dimensional 
family of SLEs (indexed by a positive real parameter $\kappa$), and they are 
the only possible candidates for the scaling limits of interfaces for 
two-dimensional critical systems that are believed to be conformally invariant. This definition of 
SLE via Loewner's equation is a dynamic one-dimensional construction: One basically describes the law of 
$\eta[t, t+dt]$ given $\eta [0,t]$, and integrates this with respect to $t$.
See e.g. \cite {Schramm1,RS,LSW1,Wstf} for an introduction to (chordal) SLE.
Furthermore, it is worthwhile stressing that this construction describes one interface, corresponding to 
specific boundary conditions in the discrete model, but that it does in general not give immediately 
access to the ``complete scaling limit'' of the system. 
This raises the following question: Is there a simple and natural way to define at once a whole 
family of SLE loops in a domain that might describe simultaneously all boundaries of clusters?

In \cite {LSWr}, a different characterization of the SLE$_{8/3}$ random curve was derived. It is shown to be the 
unique random curve in a domain, that satisfies a certain conformal restriction property. This 
characterization is ``global'' and does not use (directly) the Markovian property. It also enabled
to identify this curve with the outer boundary of a certain reflected Brownian motion \cite {LSWr},
and with the outer boundary of a certain union of Brownian excursions \cite {Wcrr}.
Hence, it is geometrically possible to construct SLE$_{8/3}$ from planar Brownian motions (recall also 
that SLE$_{8/3}$ is conjectured to be the scaling limit of the half-plane self-avoiding walk \cite {LSWsaw}).

When $\kappa =2$, SLE has been proved in \cite {LSWlesl}
to be the scaling limit of the loop-erased random walk. This gives 
a heuristic justification to the fact proved in \cite {LSWr} that adding Brownian loops 
in a proper way on the top of an SLE$_2$ curve gives the same hull as a Brownian motion.
In fact, a similar result holds for all $\kappa \in (0, 8/3)$: Adding Brownian loops to 
an SLE$_\kappa$ gives a sample of a conformal restriction measure, as defined in \cite {LSWr}.
This fact can be related to some representations of infinite-dimensional Lie Algebras \cite {FW}.
In a way, this shows that SLE$_\kappa$ 
for $\kappa \in (0,8/3]$ could also be characterized implicitly and globally via planar Brownian 
motions: It is the only simple curve such that if one adds a certain density of Brownian loops, one gets the same hull 
as the union of some Brownian motions (all these aspects relating restriction measures to SLE are 
reviewed in \cite {Wcrr}).
This raises naturally the following question: What can one say  for $\kappa \in (8/3, 4]$, which is in fact the 
physically more interesting part (in CFT language, it corresponds to positive central charge) that
is supposed to  
correspond for instance to the scaling limit of critical Potts (i.e. random cluster) models for $q \in (1,4]$?
Is there such a relation with Brownian motions, loop soups?
As we shall see, this relation is in fact in a sense richer.

\section {The Brownian loop soup percolation.}

We now use the Brownian loop soup introduced in \cite {LWls} to define conformally invariant random fractal domains:
Start with a Brownian loop soup with small intensity $c$ in a bounded open (not necessarily 
simply connected) domain $D$ ($D = \H$ is also licit). This is a countable Poissonian collection 
of Brownian loops $(l_j, j \in J)$ that stay in $D$, and is conformally invariant:  The image of this 
loop soup under a conformal mapping $\Phi$ is a loop soup with the same intensity in the domain $\Phi (D)$.
Furthermore, it satisfies restriction: If one restricts a loop soup in $D$ to those loops that stay in 
$D' \subset D$, one gets a sample of the loop soup in $D'$.
The parameter $c$ is measuring the intensity of the Poissonian procedure: For instance, the union of two 
independent loop soups of intensity $c$ is a Brownian loop soup with intensity $2c$.

Every point in $D$ is almost surely encircled by a countable number of loops in the loop soup, but there 
exist (for small $c$) many points that are ``free'' and not encircled by any loop. In fact, one can prove that 
(when $c<10$) the dimension of the set of free points is 
almost surely $2-c/5$ (while for $c>10$, no point is free), building on 
 the relation between the loop soup and the 
Brownian bubbles derived in \cite {LWls}.
This suggests that for very small $c$, there might exist whole 
paths of points that are not encircled by any loop, i.e., paths 
that do intersect no loop in the loop soup. 

We now study the set $M:= D \setminus \cup_{j \in J} l_j$.
This Cantor-like random set is the main subject of the 
present paper. It is natural to construct the loop soup clusters:
For any two loops $l$ and $l'$ in the loop soup, we say that two loops are in the same cluster if there exist a finite 
sequence of loops $l^0=l, l^1, \ldots, l^n=l'$ in the loop soup such that for all $j \le n$,
$l^j \cap l^{j-1} \not= \emptyset$.
This defines the loop soup cluster $K(l)$ as the union of all loops $l'$ that satisfy this property.
The loop soup therefore defines a countable family $(K_i, i \in I)$ of (connected) loop soup clusters.
It is possible to show that:

\begin {proposition}
\label {p1}
There exists $c_0$ such that if $c \le c_0$, then almost surely:
All loop soup clusters are at positive distance of $\partial D$ and they
 are at positive distance from each other.
For any fixed two points $a$ and $b$ on the boundary of $D$, there exists continuous paths 
 from $a$ to $b$ in $M$ (i.e. that avoids all loops).
\end {proposition}

Note that clearly, when $c$ becomes large, since almost surely, no point is free, the statement does 
not hold: There exists only one loop soup cluster, it is dense in $D$, and its complement is 
completely disconnected.
As we shall see, it will be natural to conjecture that the critical value of $c$ is $1$.

This proposition (and its proof) has similarities with multi-scale Poisson percolation models as studied in 
 Chapter 8 of \cite {MR}, and with Mandelbrot's fractal percolation model \cite {Ma,CCD}.

\section {SLE as loop soup cluster boundaries}

When $c$ is small, it is also possible to see that the ``exterior boundary'' of $K_j$ consists of 
a union of simple (disjoint) loops. For instance, for simply connected 
$D$, one can define the outer boundary $\partial_j^{out}$
of $K_j$ as the inner boundary of the connected component 
of $D \setminus K_j$ that also has $\partial D$ on its boundary.
This associates (for small $c$) to each realization of the loop soup, a countable collection of 
simple loops $(\partial_u, u \in U)$. 
Note that even though we know \cite {LSW2}
that the dimension of the outer boundary of each Brownian loop is $4/3$, the 
curves $\partial_k^{out}$ are outer boundaries of a countable union of such loops, so that their 
fractal dimension can be different.
In fact:
\begin {proposition}
\label {t3}
When $c \le c_0$, then the curves $\partial_k$ are SLE$_\kappa$ type-curves, where $\kappa \in (8/3, 4]$ and $c$ are  
related by $c= c(\kappa)= (3\kappa -8) (6- \kappa)/ (2 \kappa)$.
\end {proposition}
This shows in particular that the Hausdorff dimension of all the curves $\partial_k$ is almost 
surely $1+ \kappa/8$ (see \cite {Bef}). 
Also, since $c(\kappa) \le 1$ for all $\kappa$, one can see that 
$c_0$ in Proposition \ref{p1} can not be larger than one. It is natural to expect that this 
proposition will hold for all $c \le 1$ (i.e. that one can take $c_0=1$, and that this is 
the critical value for loop-soup percolation).

The loop soup percolation exterior boundaries therefore define at once a countable conformally invariant 
collection of 
disjoint SLE$_\kappa$-type curves in $D$.  When one perturbs the boundary of the domain $D$ and looks
how the law of this family is changed, one can use the restriction property of the Brownian loop soup
to give explicit Radon-Nikodym derivatives between the laws in different domains, in term of  the
measures on Brownian loops/bubbles. In particular, 
this shows that it behaves as expected for a conformal field theory with central charge $c$.

The statement in the proposition is a little bit vague. 
One way to make the relation between SLE$_\kappa$ and the outer boundaries of clusters precise goes as follows:
Consider a small $c$, choose $\alpha >0$ and $\kappa \in (8/3,4]$ appropriately i.e. such that 
$\alpha = (6 - \kappa)/{2 \kappa}$ and $c=c(\kappa)$.
If $D$ is a simply connected domain and $a \not= b$ two boundary points (one can take for instance
$D= \H$, $a=0$ and $b= \infty$), we define a random simple curve $\gamma$ as a sample
of the one-sided restriction measure (from $a$ to $b$ in $D$) with exponent $\alpha$. 
This random curve is defined in \cite {LSWr} and can be viewed as the boundary of 
a certain reflected Brownian motion from $a$ to $b$, or as the boundary of a union of a Poissonian sample
Brownian excursions (see \cite {Wcrr}).
Attach to $\gamma$  the union of all the loop soup clusters 
(of an independent loop soup with intensity $c$) that it intersects. 
Call the right-boundary of this set $\eta$. Then:

\begin {proposition}
$\eta$ is (exactly) an SLE$_\kappa$ curve.
\end {proposition}

In fact, if one chooses $c$ as before, but starts with another $\alpha$, one gets a so-called SLE($\kappa, \rho$)
curve, as defined in \cite {LSWr}, where
$\alpha (\kappa, \rho) = (\rho +2) (\rho +6 - \kappa) / {4 \kappa}$.

The idea of the proof goes as follows: Studying the loop soup itself (see above) shows that
$\eta$ is a simple curve that stays away from the boundary of $D$.
Because of the restriction properties of both the loop soup
and the curve $\gamma$, it is possible to argue that the curve $\eta$ satisfies a ``Markovian-type''
property that basically implies that it is an SLE($\kappa, \rho$) process for some $\kappa$ and $\rho$
(this uses the fact that the Bessel processes used to define these SLE processes are the 
only real continuous Markov process with Brownian scaling).
The values of $\kappa$ and $\rho$ can then be worked out by looking at how the law of $\eta$ behaves when 
one perturbs the boundary of the domain $D$, and comparing this with the local martingales pointed out 
for SLE($\kappa,\rho$) by Dub\'edat in \cite{Dub}.

One consequence of this result is the ``reversibility'' of these SLE, i.e.: $\eta$ and $-1/\eta$ have the 
same law (if one looks at SLE from the origin to infinity in $\H$).

Along similar lines, it is possible to give a heuristic justification to the fact that
chordal SLE$_\kappa$ is the (annealed) uniform measure on self-avoiding walks on the 
random fractal $M$, when $c=c(\kappa)$:  
Consider a Brownian loop soup
in $D$ with 
intensity $c$. Suppose that 
there exist a conformally invariant and uniform way to choose a simple curve
$\eta$ from $a$ to $b$ that avoids all loops in the loop soup. 
This is  very vague and as unprecise as to say that there exists a uniform
conformally invariant way to choose a self-avoiding curve from the origin to infinity, which would 
correspond to the case $c=0$; see \cite {LSWsaw,Wcrr} for how one can (and cannot) make such definitions
rigorous. In fact, one should rather speak of ``intrinsic'' measure rather than uniform measures 
(see \cite {Wcrr}). Anyway, if such a definition would hold, then 
the conformal restriction property of the loop soup implies readily that this curve satisfies 
the Markovian property, so that it should be an SLE curve. The value $\kappa$ can then be 
determined as before studying the way in which the law of this curve is changed under restriction. 
It is worthwhile exploring whether 
this is related to the ``quantum gravity'' approach to critical phenomena developed by physicists,
see e.g. \cite {KPZ} (here, we interpret SLE as self-avoiding walk on a natural continuous random geometric object).

\section {Consequences}
This has many consequences. We plan to address the following items in forthcoming papers:
\begin {itemize}
\item
The restriction property of the loop soup shows that these clusters satisfy a ``Markovian property'' in space.
This is related to various things, in particular to certain representations of the Virasoro 
algebra.
\item
This family of SLE loops defines at once a big family of observables. This allows to define a natural 
$L^2$ space, on which the Loewner semi-group acts. This is related to considerations from conformal 
field theory, as for instance in \cite {BPZ,Ca}. 
\item
One can do similar things using radial restriction measures and radial restriction.
\item
This construction works obviously on any Riemann surface, and the value of the critical $c$ is the 
same as on the plane.
This is of course related to some of the previous items. 
\end {itemize}
We now conclude the paper with some comments:
\begin {itemize}
\item
There exists a representation of correlation functions for spin systems (see e.g. \cite {BFS}) via random walks.
Maybe  this is related to the present loop soup percolation (i.e. chain of Brownian loops) representation, and can have fruitful consequences.
\item
The $c=0$ model (i.e. the scaling limit of percolation clusters, say)
does not appear easily here. This is probably related to the fact that in the CFT approach, the 
$c=0$ case is often treated via the $c \to 0$ limit \cite {Ca}.
\item
The construction of SLE in proposition \ref {t3}
is very non-symmetric, and one may be surprised to obtain a symmetric curve in the end (with respect to the imaginary axis). Note that it the limiting case $c=0$ (where no loop is present), this was already proved 
to be the case \cite {LSWr}.
\item
The construction of SLE in proposition \ref {t3} is very ``two-dimensional'' and in spirit very 
different from the Loewner equation approach. The fact that these two 
constructions are equivalent is not so surprising after all: Because of 
the fact that there are only few candidates to describe conformally invariant models in two-dimensions, 
many a priori different definitions in fact coincide (which is one of the 
instrumental observations \cite {LW2} that allowed to use SLE to determine the Brownian exponent \cite{LSW1,LSW2}). As the Brownian loop soup is a rich conformally invariant object, it is not so surprising that it ``contains'' SLE curves. The same remark 
applies to the relation between the Gaussian Free Field and SLE curves/loops recently
discovered by Oded Schramm and Scott Sheffield \cite {SS}. It also raises the question of whether there exists 
a direct link between the Gaussian Free Field and the loop soup.
\end {itemize}
 
\section* {Acknowledgements}
I thank Vincent Beffara, Greg Lawler and Nikolai Makarov for stimulating and useful discussions. I also acknowledge the 
hospitality of the Isaac Newton Institute, where part of this work was done.

\end {document}